
%
%
%
%
%
%
\input amstex
\documentstyle{amsppt}
\magnification=1200
\tolerance=10000
%
%
%
%
%
%
\def\R{\Bbb R}
\def\Z{\Bbb Z}
\def\N{\Bbb N}

\def\s{\sigma}
\def\w{\omega}
\def\tens{\widehat \otimes}

\def\te{\tilde e}
\def\tp{\tilde \pi}
\def\pa{\parallel}
\def\vr{\vec r}
\def\Finf{F^{\infty}}
%
%
\def\schwr{\Cal S (\R)}
%
%
%

\def\GLA{G \rtimes^{\w} A}
\def\GS{G\rtimes \schwm}
%
%
\def\cc{C_{c}^{\infty}(G)}
%
%
%

\def\schwm{\Cal S (M)}
\def\schwmh{\Cal S_{H}^{\sigma}(M)}
%
%
%
%
%
%
%
%
%
%
\hoffset=0.2truein
\voffset=0.2truein
\hsize=7.2truein
\vsize=8.7truein
\addto\tenpoint{\normalbaselineskip=20pt\normalbaselines}
\addto\eightpoint{\normalbaselineskip=16pt\normalbaselines}
\topmatter
\title
A Factorization Theorem for Smooth Crossed Products
\endtitle
\author
Larry B. Schweitzer
\endauthor
\address
Department of Mathematics and Statistics, University of Victoria,
Victoria, B.C., Canada V8W 3P4  \endaddress
\email lschweit\@sol.uvic.ca \endemail
\endtopmatter
\document
\heading Introduction \endheading
\par
By a remarkable theorem of Dixmier and Malliavin \cite{DM, Theorem 3.3},
it is known that the convolution algebra $C_{c}^{\infty}(G)$
of compactly supported $C^{\infty}$-functions on a Lie group $G$ satisfies
the factorization property, namely that every set of $C^{\infty}$-vectors
$E$ for the action of $G$ is equal to the finite
linear span $C_{c}^{\infty}(G)E$.
In this paper, we replace $C_{c}^{\infty}(G)$ by the smooth
crossed products for transformation groups $G\rtimes \schwm$
defined in \cite{Sc 1}.
We define an appropriate notion of a differentiable $\GS$-module,
which generalizes the notion of $C^{\infty}$-vectors for actions of
Lie groups. (This definition was first introduced by F. Du Cloux \cite{DuC 1}
\cite{DuC 2}.)
Under the assumption that the Schwartz functions $\schwm$ vanish
rapidly with respect to a continuous,
proper map $\s \colon M \rightarrow [0, \infty)$, we then
show that $\GS$ satisfies the factorization property,
namely that any differentiable $\GS$-module $E$ is the finite
span of elements of the form $ae$ where $ a \in \GS$  and $e \in E$.
In the course of doing this, we also show that if a Fr\'echet algebra
$A$ has the factorization property, then the smooth crossed
product $G\rtimes A$ does also.
\par
Other aspects of the  representation theory
of the smooth crossed products $\GS$ are studied in \cite{Sc 2}.
I would like to thank Berndt Brenken for a pleasant stay at the University of
Calgary, where I wrote the first draft of this paper.
\heading \S 1 Differentiable Representations and Multipliers \endheading
\par
We define what it means for an algebra and a representation to
be differentiable.
We shall use representation and module
terminology interchangeably throughout this paper.
Everything we do will be for left modules, though
similar statements are also true for right modules.
\subheading{Definition 1.1}  By a {\it Fr\'echet algebra} we mean
a Fr\'echet space with an algebra structure for which the
multiplication is jointly continuous.
(We do not assume Fr\'echet algebras are $m$-convex.)
Let $A$ be a Fr\'echet algebra.
By a Fr\'echet $A$-module, we mean a Fr\'echet space
$E$ that is an $A$-module for which the map $(a, e)\mapsto ae$
is jointly continuous.
The $A$-module
$E$ is {\it non-degenerate}
({\it differentiable}) if $\bigl\{ v\in E \,\bigl|\,
Av = 0 \bigr\} =  \bigl\{0\bigr\}$, and
the image of the canonical map $A {\widehat \otimes} E \longrightarrow E$
is dense (onto) \cite{DuC 2,  D\'efinition 2.3.1}.
(All tensor products will be completed in the projective
topology.)
We make the same definitions for right Fr\'echet $A$-modules.
If $A$ is differentiable both as a
left and right $A$-module, then we say that the Fr\'echet algebra $A$ is
{\it self-differentiable }.
\par
(In \cite{DuC 2, D\'efinition 2.3.1},  a self-differentiable
Fr\'echet algebra is called a
\lq\lq differentiable Fr\'echet algebra\rq\rq.  However, this terminology
would suggest that the algebra has a derivation acting
on it, or that it is a set of $C^{\infty}$-vectors for
the action of a Lie group.    This is not the case;
any C*-algebra is a \lq\lq differentiable Fr\'echet algebra\rq\rq\
since any element of the algebra can be written as a linear combination
of four positive elements, each of which has a square root.
So I prefer to say \lq\lq self-differentiable\rq\rq\
instead of \lq\lq differentiable\rq\rq.)
\par
If $E$ is non-degenerate, we let $E_{s}(A)$ be the
image of the canonical map $A\widehat \otimes E \longrightarrow E$.
Then $E_{s}(A)$ inherits the quotient topology from
$A\widehat \otimes E$, making
$E_{s}(A)$ a Fr\'echet $A$-module.
When $A$ is self-differentiable, the
$A$-module $E_{s}(A)$ is always differentiable \cite{DuC 2, Lemme 2.3.4}.
\par
If $G$ is a topological group, we say that a
Fr\'echet space
$E$ is
a {\it continuous $G$-module} if $G$
acts on $E$ by
continuous automorphisms, and
for each $e \in E$ and each continuous seminorm $\pa \quad \pa$
on $E$, the map $g \mapsto \pa ge \pa$ is continuous.
If $G$ is a Lie group, we  say that $E$ is a {\it differentiable}
$G$-module if the action of $G$ on $E$ is differentiable.
\par
We say that a self-differentiable Fr\'echet algebra $A$
satisfies the {\it factorization property} if every differentiable
$A$-module $E$ is the finite span of elements of the form $ae$
where $a \in A$ and $e \in E$.
Note that in particular $A$ will be the finite span of products
of elements of $A$.
\par
Note that if $A$ is a unital Fr\'echet algebra (this corresponds to the
case of the group algebra of a discrete group),
then $A$ is self-differentiable, every
$A$-module is differentiable, and $A$ satisfies the factorization
property.
\par
If $G$ is a compact Lie group, and $A= C^{\infty}(G)$
is the convolution algebra of $C^{\infty}$ functions
on $G$, then an $A$-module $E$ is differentiable
if and only if the action of $G$ on $E$ is differentiable
(see Theorem 5.3 below or \cite{DuC 2, Exemple 2.3.3}).
It follows immediately from \cite{DM, Theorem 3.3}
that Schwartz functions  $\schwr$ on $\R$  with convolution multiplication
is a self-differentiable Fr\'echet algebra,
since the canonical map $\schwr \otimes \schwr \longrightarrow \schwr$
is onto. (Here $\otimes$ denotes the algebraic tensor product.)
\subheading{Example 1.2}  We give an example of a self-differentiable
Fr\'echet (in fact Banach) algebra without the factorization
property.  Let $A=l_{1}(\Z)$ with pointwise multiplication.
Then $c_{f}(\Z)$  is dense in $A$ so $A$ is nondegenerate.
Since $A{\widehat \otimes} A \cong l_{1}(\Z \times \Z)$,
and the canonical map $\pi \colon A{\widehat\otimes}A \longrightarrow
A$ is given by evaluation along the diagonal, $A$ is self-differentiable.
\par
A quick calculation shows that $\pa \varphi * \psi \pa_{1/2} \leq
\pa \varphi \pa_{1} \,\pa \psi \pa_{1}$ and $\pa \varphi + \psi \pa_{1/2}
\leq 2\bigl( \pa \varphi \pa_{1/2} + \pa \psi \pa_{1/2} \bigr)$.
Hence the algebraic span $A^{2}$ is contained in $l_{1/2}(\Z)$.
Since $l_{1/2}(\Z) \not= A$, the algebra $A$ does not have the
factorization property.
\par
Similar arguments show that the Banach algebra $l_{2}(\Z)$ with
pointwise multiplication is an example of a  nondegenerate
but non self-differentiable Banach algebra.
\subheading{Definition 1.3} We say that $T$ is a {\it multiplier}
for a Fr\'echet algebra $A$ if $T$ acts as a continuous linear operator
both on the left of $A $ and the right of $A$, and
the left and right actions commute.  It follows that
 for every seminorm $\pa \quad \pa_{d}$ on $A$, there is some
$C>0$ and another seminorm $\pa \quad \pa_{m}$ on $A$ such that
$$ \max\biggl( \pa Ta \pa_{d},\,
\pa aT\pa_{d} \biggr)
 \leq C \pa a \pa_{m}, \qquad a \in A.  $$
\subheading{Example 1.4}  In general if $T$ is a multiplier,
and $E$ is a nondegenerate $A$-module, the action of $T$ on $A$
does not extend to an action on $E$. For example, let $A$ be
Schwartz functions $\schwr$ on $\R$ with pointwise multiplication,
and let $T$ be multiplication by the function $r^{2}$.
Let $E$ be the nondegenerate $A$-module $C_{0}(\R)$ with action of
$A$ given by pointwise multiplication. (Here $C_{0}(\R)$
denotes the set of continuous functions on $\R$ which vanish
at infinity).  Let $f$ be any continuous
function which vanishes like $1/r^{2}$
at infinity on $\R$.  Then $f \in E$ and $Tf$ does not vanish at
infinity so $Tf \notin E$.
\proclaim{Theorem 1.5} Let $A$ be a self-differentiable Fr\'echet algebra
and let $E$ be a differentiable  $A$-module.
Let $T$ be a  multiplier for $A$.  Then there is a unique
action of $T$ on $E$ as a continuous linear operator, which is consistent
with the action of $A$ on $E$. \endproclaim
\demo{Proof}
Since $T$ is a continuous linear map from $A$ to $A$, $T$ also gives
a continuous linear map of the projective completions
$T\colon A{\widehat \otimes} E \longrightarrow A{\widehat \otimes} E$
\cite{Tr, Proposition 43.6}.
Since $E$ is differentiable, to see that this map induces a continuous
linear map on $E$, it suffices to show that $T$ leaves the kernel of
the canonical map $\pi\colon A{\widehat \otimes} E \longrightarrow E$
invariant.
Assume that $\pi(\eta)=0$ for $\eta \in A{\widehat \otimes} E$.
Let $b \in A$.  Then
$$b \pi(T\eta)= \pi( bT\eta) = bT\pi(\eta) = 0,$$
since $bT \in A$.  Hence $\pi(T\eta)\in E$ is annihilated by every
element of $A$. Since $E$ is a nondegenerate $A$-module,
it follows that $\pi(T\eta)=0$.  Hence $T$ leaves the kernel of
$\pi$ invariant.
\qed
\enddemo
\heading \S 2 Smooth Crossed Products \endheading
\par
We recall the definitions of our smooth crossed products from \cite{Sc 1}.
First, let $H$ be a Lie group and let $M$ be  a locally compact space
on which $H$ acts.  We say that a Borel measurable function
$\s \colon M \rightarrow [0, \infty)$ is a {\it scale} if it is
bounded on compact subsets of $M$.
We say that a scale $\s$ {\it dominates} another scale $\gamma$
if there exists $C, D>0$ and $d\in \N$ such that
$ \gamma(m) \leq C\s(m)^{d} + D$ for  $ m \in M$.
We say that $\s$ and $\gamma$ are {\it equivalent
($\s \thicksim \gamma$)} if they dominate each
other.  We may always replace a scale $\s$ with the equivalent
scale $1+ \s$, so that we lose no generality by assuming $\s \geq 1$.
{}From now on, we will assume this.
If $h \in H$, define
$\s_{h}(m) = \s(h^{-1}m)$.   We say that $\s$ is {\it uniformly
$H$-translationally equivalent} if
for every compact subset $K$ of $H$ there exists $C_{K}>0$ and $d \in\N$
such that
$$ \s_{h} (m) \leq C_{K}\s(m)^{d},\qquad m \in M,\, h \in K. \tag 2.1 $$
If $\s$ is a uniformly $H$-translationally equivalent scale on $M$,
we may define the {\it $H$-differentiable $\s$-rapidly vanishing
functions $\schwmh$} by
$$ \schwmh = \biggl\{ f \in C_{0}(M), \,
\text{$f$ $H$-differentiable} \, \biggr| \, \pa \s^{d}X^{\gamma}f \pa_{\infty}
<\infty \text{ and } X^{\gamma}f \in C_{0}(M) \,
\biggr\}, $$
where $X^{\gamma}$ ranges over all
differential operators from the Lie algebra of
$H$, and $d$ ranges over all natural numbers.
We topologize $\schwmh$ by the seminorms
$$ \pa f \pa_{d, \gamma} = \pa \s^{d}X^{\gamma}f \pa_{\infty}. $$
Then $\schwmh$ is a Fr\'echet *-algebra
under pointwise multiplication,
with differentiable action of $H$ \cite{Sc 1, \S 5}.
\par
Next, let $G\subseteq H$ be a Lie group with differentiable
inclusion map $\iota \colon G \hookrightarrow H$.  Let $\w\geq 1$ be
a scale on $G$.
Let $E$ be any Fr\'echet space.
 We define the {\it differentiable $\w$-rapidly
vanishing functions $\Cal S^{\w}(G, E)$ from $G$ to $E$} to be the set
of differentiable functions $\varphi$
from $G$ to $E$ such that
$$\pa \varphi \pa_{d, \gamma, m} =
\int_{G}\pa\w^{d}X^{\gamma}\varphi(g) \pa_{m} dg <\infty, \tag 2.2$$
where $X^{\gamma}$ is any differential operator from the Lie
algebra of $G$ acting by left translation,
$\pa \quad \pa_{m}$ is any seminorm for $E$,
and $d$ is any natural number.
We topologize $\Cal S^{\w}(G, E)$ by the seminorms (2.2).
\par
We say that the action of $G$ on a $G$-module $E$ is {\it $\w$-tempered}
if for every $m \in \N$ there exists $C>0$, $d \in \N$
and $k \in \N$ such that
$$ \pa \alpha_{g}(e)
\pa_{m} \leq C\w(g)^{d} \pa e \pa_{k}, \qquad e \in E, \,g \in G. $$
Simple arguments show that every closed $G$-submodule and
every quotient of a tempered $G$-module is again a tempered
$G$-module.
We say that $\w$ is {\it sub-polynomial} if there exists $C>0$,
$d\in \N$ such that
$$ \w(gh) \leq C\w(g)^{d} \w(h)^{d}, \qquad g,\, h \in G. $$
The {\it inverse scale} $\w_{-}$ is defined by $\w_{-}(g) = \w(g^{-1})$.
We say that $\w$ {\it bounds $Ad$ on $H$} if there exists $C>0$, $d\in \N$
such that
$$ \pa Ad_{g} \pa \leq C\w(g)^{d}, \qquad g \in G,  $$
where $\pa Ad_{g} \pa$ is the operator norm of $Ad_{g}$ as an operator
on the Lie algebra of $H$.
And finally, if $\w$ is a sub-polynomial scale on $G$ such that
$\w_{-}$ bounds $Ad$ on $H$, and $\s$ satisfies
$$ \s(gm) \leq C\w(g)^{d}\s(m)^{l}, \qquad g\in G,\quad  m\in M\tag 2.3$$
for some $C>0$ and $d, l \in \N$, then
we say that {\it $(M,\s, H)$ is a
scaled $(G, \w)$-space}.
\proclaim{Theorem 2.4 \cite{Sc 1, Theorem 2.2.6,  Theorem 5.17}}
Let $\w$ be a sub-polynomial
scale on a Lie group $G$ such that
$\w_{-}$ bounds $Ad$ on $G$.
Assume that the action of $G$ on a Fr\'echet algebra $A$ is
continuous and $\w$-tempered.
Then $\Cal S^{\w}(G, A)$ is a Fr\'echet algebra under convolution,
which we denote by $G\rtimes^{\w} A$.
\par
Moreover, if $(M, \s, H)$ is a scaled $(G, \w)$-space,
then the action of $G$ on $\schwmh$ is differentiable and
$\w$-tempered.  In particluar, $G \rtimes^{\w} \schwmh$ is a Fr\'echet
algebra under convolution.
\endproclaim
\par
See \cite{Sc 1, \S 5} or \cite{Sc 2} for examples.
\vskip\baselineskip
\heading \S 3  Differentiable Scales for $M$
\endheading
\proclaim{Theorem 3.1} Every  uniformly $H$-translationally equivalent
scale $\s$ on $M$ is equivalent to an $H$-differentiable scale $\tilde \s$
on $M$ for which there is some $d \in \N$ such that
for every differential operator $X^{\gamma}$ from the Lie algebra
of $H$ we have
$$(\exists C_{\gamma}>0) \quad
X^{\gamma}{\tilde \s}(m) \leq C_{\gamma}
{\tilde \s}^{d}(m), \qquad m \in M. \tag 3.2$$
If $\s$ is continuous to begin with, then
the scale $\tilde \s$ produced in the
proof is also.
\endproclaim
\demo{Proof} Let $\s$ be any  uniformly
$H$-translationally equivalent scale on $M$.
Let $K$ be a compact neighborhood of $e$ in $H$ such that
$$ \s_{g}(m) \leq C\s(m)^{d}, \qquad m \in M \tag 3.3 $$
and
$$ \s(m) \leq C\s_{g}(m)^{d} \qquad m \in M $$
for every $g \in K$.
Let $\varphi \in C_{c}^{\infty}(H)$ be any nonnegative function
with support contained in $K$ such that $\int \varphi(g) dg =1$.
Define
$$ {\tilde \s}(m) = \int \varphi(g) \s_{g}(m) dg. $$
Then ${\tilde \s}(m) \geq 1$ and $\tilde \s$ is Borel measurable on $M$.
If $\s$ is continuous, then taking limits inside the integral
shows that $\tilde \s$ is also.
We show that $\s \thicksim {\tilde \s}$.
By (3.3), we have
$$ {\tilde \s}(m) \leq \int \varphi(g) C\s(m)^{d} dg = C\s(m)^{d}  $$
so $\s$ dominates $\tilde \s$ (in particular, $\tilde \s$
is bounded on compact sets).  Similarly,
$$ \s(m)^{1/d} = \int \varphi(g) \s(m)^{1/d} dg
\leq \int \varphi(g) C^{1/d}\s_{g}(m)dg
= C^{1/d} {\tilde \s}(m),\tag 3.4 $$
so
$\s(m) \leq C {\tilde \s}^{d}(m)$.
\par
We show that $\tilde \s$ is differentiable, and that the derivatives
satisfy the bounds (3.2).
We have
$$ X^{\gamma}{\tilde \s}(m) = \int X^{\gamma}\varphi(g) {\s}_{g}(m) dg.$$
So ${\tilde \s}(m)$ is an $H$-differentiable function on $M$.
Using (3.3), we bound the derivative
$$ |X^{\gamma}{\tilde \s}(m)|
\leq \int \, |X^{\gamma}\varphi(g)| \, C\s(m)^{d} dg
= C_{\gamma} C \s(m)^{d}. $$
Since $\tilde \s$ dominates $\s$ (see (3.4)), we have (3.2)
\qed
\enddemo
\par
We say that a scale $\s \colon M \rightarrow [0, \infty)$ is {\it proper}
if the inverse
image $\s^{-1}(K)$ of every compact subset $K$ of $[0, \infty)$
is relatively compact.  The property of being proper is preserved under
equivalence.
\proclaim{Proposition 3.5}
Let $\s$ be a continuous uniformly $H$-translationally equivalent
$H$-differentiable
scale on $M$ with property (3.2).
Then  $\s$ is
a multiplier on $\schwmh$.  If $\s$ is proper,
then there is a natural continuous algebra homomorphism
$\schwr\rightarrow \schwmh$
given by $\varphi \mapsto \varphi\circ \s$.
\endproclaim
\demo{Proof}
To see that $\s$ is a multiplier on $\schwmh$, let $f \in \schwmh$.
Then
$$\pa \s^{l}X(\s f)\pa_{\infty}=
\pa\s^{l}((X\s)f + \s(Xf))\pa_{\infty}
\leq \pa \s^{l}C\s^{d} f \pa_{\infty} + \pa \s^{l+1} Xf \pa_{\infty}.$$
Similar arguments show that for higher derivatives, we also have
$\pa \s^{l} X^{\gamma}(\s f) \pa_{\infty}$ bounded by some linear
combination of seminorms of $f$.
The function $X^{\gamma}(\s f)$ is a continuous function on $M$,
since $\s$ and $f$ are continuous and $H$-differentiable.
Since for each $l \in \N$ and $\gamma$, the function
$\s^{l} X^{\gamma} f$ vanishes
at infinity \cite{Sc 2, Proof of Proposition 5.2},
and $|X^{\gamma} \s | \leq C_{\gamma} \s^{d}$
by (3.2), we also have $X^{\gamma} (\s f) \in C_{0}(M)$
for all $\gamma$.
Hence $\s f \in \schwmh$ and $\s$
is a multiplier on $\schwmh$.
\par
For the second statement,
it suffices to show that the seminorms of $\varphi \circ \s$ in
$\schwmh$
are  bounded by linear combinations  of
seminorms of $\varphi $ in $\schwr$,
and also that $X^{\gamma}(\varphi \circ \s) \in C_{0}(M)$.
We apply the chain rule.  If $X$ is in the Lie algebra of
$H$, then
$$ X (\varphi \circ \s )(m) = (\varphi^{\prime}\circ \s)(m)
X\s(m).  $$
So
$$
\pa \s^{l}X\varphi\circ \s \pa_{\infty}
=\sup_{m\in M} |\s(m)^{l} X(\varphi \circ \s)(m)|\leq
\sup_{m\in M}|\s^{l}\varphi^{\prime}(\s(m))C\s(m)^{d}|
\leq \sup_{r\in \R}|r^{l}\varphi^{\prime}(r)Cr^{d}|,
$$
where the last expression is a seminorm of $\varphi$
in $\schwr$.  Since $\s$ is proper, and $ r^{d}\varphi^{\prime}(r)$
vanishes at infinity,
the fact that
$|X(\varphi \circ \s)(m)| = |\varphi^{\prime} (\s(m)) X\s(m)|
\leq |C \s(m)^{d} \varphi^{\prime}(\s(m))|$ implies
$X(\varphi \circ \s) \in C_{0}(M)$.
Similar arguments work for higher derivatives.
\qed
\enddemo
\subheading{Example 3.6}
We consider the case when $M= \R^{n}$, and $\s(\vr) =
r_{1}^{2} + \dots r_{n}^{2}$.
Then $\s$ is a differentiable scale on $\Cal S(\R^{n})$.  The map
$\schwr \rightarrow \Cal S(\R^{n})$ in Proposition 3.5 above is given by
$\varphi(\vr) = \varphi(r_{1}^{2}+ \dots r_{n}^{2})$.  The image of this
map consists of radially symmetric functions on $\R^{n}$.  In this sense,
a differentiable scale can be regarded as a generalized
\lq \lq radial\rq \rq\
coordinate for $M$, and the image of $\schwr$ in $\schwmh$
consists of functions
which depend only on this radial coordinate.
\heading \S 4 Factorization Property for
$\schwmh$ \endheading
\par
We recall some of the functions on $\R$ defined in the proof of the
Dixmier-Malliavin theorem \cite{DM}.
Let $\lambda= (\lambda_{0}, \lambda_{1}, \dots, \lambda_{k}, \dots )$
be any subsequence of $(1, 2,\dots, 2^{k}, \dots)$.
For $x\in\R$, let
$$\varphi_{\lambda}(x) = \prod_{k=0}^{\infty}
\biggl( 1 + {x^{2}\over \lambda_{k}^{2}} \biggr), \qquad
\chi_{\lambda}(x) = \varphi_{\lambda}(x)^{-1}.$$
We show that $\varphi_{\lambda}$ is a well defined function from
$\R $ to $[1, \infty)$.
For $k$ sufficiently large we have $x^{2}<\lambda_{k}^{2}$, so
$$ 1+ {x^{2}\over {\lambda_{k+p}^{2}}} \,<\, 1+{1\over{2^{p}}} $$
for all $p \in \N$.  Since
$$ \log\biggl(\,\biggl( 1+ {1\over{1^{}}}\biggr)
\biggl( 1+ {1\over{2^{}}}\biggr)
\dots
\biggl( 1+ {1\over{2^{p}}}\biggr)
\dots\,\biggr) =
\sum_{p=0}^{\infty} \log\biggl(1+{1\over{2^{p}}}\biggr)
\leq \sum_{p=0}^{\infty}{1\over{2^{p}}} <
\infty, $$
we see that $\varphi_{\lambda}(x)$ is well defined for any $x \in \R$.
%
\par
It is shown in \cite{DM, \S 2.3} that $\chi_{\lambda}$ is in $\schwr$.
Also, it is shown in the proof of \cite{DM, Lemme 2.5 - p. 309-310}
that for any sequence $(\beta_{0}, \beta_{1}, \dots)$ of positive numbers,
there exists a sequence $(\alpha_{0},\alpha_{1}, \dots)$
of positive numbers,
and a sequence $(\lambda_{0}, \lambda_{1}, \dots)$ as above
such that  $\alpha_{n}$ occur in the expansion
$$\varphi_{\lambda}(x)
= \sum_{n=0}^{\infty} \alpha_{n} x^{2n} $$ and satisfy
$$ \alpha_{n} \leq \min(\beta_{n}, 1/n^{2}). \tag 4.1$$
We use this to show that $\schwmh$ has the factorization property.
\proclaim{Theorem 4.2}  Assume that $\s$ is continuous and proper.
Then for every function $\psi \in \schwmh$, there are
$\theta, \phi \in \schwmh$ such that $\psi = \theta\phi$.
\endproclaim
\demo{Proof}
Let $\psi \in \schwmh$, and let $\s$ be an $H$-differentiable
scale as in Theorem 3.1 above. Define
$$ M_{d,l, n} = \max_{|\gamma|\leq l } \pa
\s^{(d+1)l}
\s^{2n}X^{\gamma}\psi \pa_{\infty}.$$
Choose $\lambda=(\lambda_{0}, \dots \lambda_{1}, \dots )$
so that the sequence $(\alpha_{0}, \alpha_{1},\dots)$ satisfies
$$ \sum_{n=0}^{\infty} \alpha_{n}M_{d,l,n} \,< \,\infty,
\qquad d,l \in \N.\tag
4.3$$
Recall that $\sum_{n=0}^{\infty} \alpha_{n} x^{2n}$ is the expansion for
$\varphi_{\lambda}$. Define ${\tilde \varphi}_{\lambda} = \varphi_{\lambda}
\circ \s\colon M \rightarrow \R$.
Then  the series for ${\tilde \varphi}_{\lambda} \psi$
converges absolutely in $\schwmh$ to an element of $\schwmh$.
For if $|\gamma|\leq l$, we have
$$ \aligned \pa &\s^{d} X^{\gamma}
\biggl( \sum_{n=k}^{\infty} \alpha_{n} \s^{2n}\biggr)
\psi \pa_{\infty}\leq \sum_{n=k}^{\infty}\alpha_{n} \pa \s^{d}X^{\gamma}
(\s^{2n}\psi) \pa_{\infty}\\
& \leq \sum_{n=k}^{\infty} \alpha_{n}\sum_{|\beta_{1}|+\dots
|\beta_{2n+1}|\leq l}
D\pa \s^{d}
X^{\beta_{1}}\s \dots X^{\beta_{2n}}\s X^{\beta_{2n+1}}
\psi \pa_{\infty} \\
& \leq \sum_{n=k}^{\infty} \alpha_{n} DC\max_{|\beta|\leq l}
\pa \s^{d}
\s^{dl}\s^{2n-l} X^{\beta}\psi \pa_{\infty}\qquad {\text{by (3.2)}} \\
&\leq \sum_{n=k}^{\infty} \alpha_{n} DC
M_{d,l,n}\qquad {\text{since $\s^{2n-l}\leq \s^{2n}$.}}\endaligned
\tag 4.4 $$
By our constraint on the $\alpha_{n}$'s (4.3), the last sum tends to
zero as $k\rightarrow \infty$.  Hence ${\tilde \varphi}_{\lambda}\psi$
converges to some well defined element $\phi\in \schwmh$.   Let $\theta
= \chi_{\lambda}\circ \s$.  By Proposition 3.5 and since $\chi_{\lambda}
\in \schwr$, we know $\theta \in \schwmh$.
Since $\chi_{\lambda}(x) = \varphi_{\lambda}^{-1}(x)$,
we have $\theta (m) {\tilde \varphi}_{\lambda}(m) = 1$ for each $m\in M$,
where $1$ denotes the identity multiplier on $\schwmh$.
It follows that $\theta(m) \phi(m) = 1\psi(m)$ and the theorem is proved.
\qed
\enddemo
\par
The following corollary was part of the motivation
for this paper.
\proclaim{Corollary 4.5} If $\s$ is continuous and proper,
then the Fr\'echet algebra $\schwmh$ is
a self-differentiable Fr\'echet algebra.
\endproclaim
\demo{Proof}
This follows directly from  Theorem 4.2.
\qed
\enddemo
\proclaim{Theorem 4.6} Let $E$ be any differentiable representation
of $\schwmh$, and assume that $\s$ is continuous and proper.
Then for  every $e\in  E$, there are
$\theta\in \schwmh$ and $ f\in E$ such that $e = \theta f$.
In particular, $\schwmh$ satisfies the factorization property.
\endproclaim
\demo{Proof} We proceed very much like the proof of Theorem 4.2 above.
Let $e\in E$, and let $\s$ be an $H$-differentiable
scale as in Theorem 3.1 above. Since $E$ is a differentiable
$\schwmh$-module, $\s^{2n}e $ is a well defined element of $E$
for each $n$ (see Proposition 3.5 and Theorem 1.5).  Define
$$ M_{m, n} =  \pa
\s^{2n}e \pa_{m},$$
where $\pa \quad \pa_{m} $ is an increasing family of seminorms for $E$.
Choose $\lambda=(\lambda_{0}, \lambda_{1}, \dots )$
so that  the sequence $(\alpha_{0}, \alpha_{1},\dots)$ satisfies
$$ \sum_{n=0}^{\infty} \alpha_{n}M_{m,n} \,< \,\infty,
\qquad m \in \N.\tag 4.7$$
Recall that $\sum_{n=0}^{\infty} \alpha_{n} x^{2n}$ is the expansion for
$\varphi_{\lambda}$. Define ${\tilde \varphi}_{\lambda} = \varphi_{\lambda}
\circ \s$.
Then  the series for ${\tilde \varphi}_{\lambda} e$
converges absolutely in $E$ to an element of $E$.
For if $m \in \N$, we have
$$ \aligned \pa &
\biggl( \sum_{n=k}^{\infty} \alpha_{n} \s^{2n}\biggr)
e \pa_{m}\leq \sum_{n=k}^{\infty}\alpha_{n} \pa
\s^{2n}e \pa_{m}\\
&\leq \sum_{n=k}^{\infty} \alpha_{n}
M_{m,n}\endaligned
\tag 4.8 $$
By our constraint on the $\alpha_{n}$'s (4.7), the last sum tends to
zero as $k\rightarrow \infty$.  Hence ${\tilde \varphi}_{\lambda}e$
converges to some well defined element $f\in E$.
The remainder of the proof is just as in Theorem 4.2.
\qed
\enddemo
\heading \S 5 Factorization Property for the Crossed Product \endheading
\subheading{Definition 5.1}
Let $\w$ be a scale on a Lie group $G$, and
let $A$ be a  Fr\'echet
algebra on which $G$ acts by algebra automorphisms $\alpha_{g}$.
Assume that we have representations of $G$ and $A$ on
a Fr\'echet space $E$ such that the action of $G$ on $E$
is $\w$-tempered and differentiable,
the action of $A$ is
differentiable,  and
the covariance condition
$$ g(ae) = \alpha_{g}(a)ge, \qquad g\in G, \quad a\in A,
\quad e\in E \tag 5.2 $$
is satisfied.
We call such a representation an  {\it $\w$-tempered
differentiable covariant representation
of ($G$, $A$)}.
\proclaim {Theorem 5.3}
Let $G$ be a Lie group with sub-polynomial scale
$\w$  such that $\w_{-}$
bounds $Ad$ on $G$, and let
$A$ be a self-differentiable Fr\'echet algebra, with an
$\w$-tempered, differentiable action $\alpha_{g}$ of $G$.
Assume that we have an
$\w$-tempered differentiable covariant representation of ($G$, $A$)
on a Fr\'echet space $E$.
Then we may integrate this representation
to get a differentiable representation of  the smooth crossed product
$\GLA$ on $E$.
\par
Conversely,  if we have a differentiable representation
of $\GLA$ on $E$, then there is
an $\w$-tempered
differentiable covariant representation of ($G$, $A$) on $E$,
whose integrated form gives back the original action
of $\GLA$ on $E$.
\par
It follows from the proof that
the smooth crossed product $\GLA$ is a self-differentiable
Fr\'echet algebra.
\endproclaim
\demo{Proof}
To simplify notation, we let $B=\GLA$.  We define an action of the
algebra $B$ on $E$ by
$$ Fe = \int_{G} F(g) (ge) dg. \tag 5.4 $$
We estimate
$$ \aligned \pa Fe \pa_{d} & \leq \int_{G} \pa F(g)ge \pa_{d}dg\\
&\leq \int_{G} \pa F(g) \pa_{m} \pa ge \pa_{k} dg \qquad
{\text{$E$ cont $A$-module}}\\
&\leq C\int_{G} \pa F(g) \pa_{m} \w(g)^{l}\pa e \pa_{n} dg
\qquad {\text{$G$ acts tempered}} \\
&\leq C  \pa F \pa_{l, m} \pa ge \pa_{n} \endaligned
\tag 5.5
$$
where $\pa F \pa_{l, m}$ are  seminorms for $B$.
So
(5.4) is well defined and continuous.
By the covariance condition (5.2), it follows that
$(F_{1}*F_{2})e=F_{1}(F_{2}e)$,
so
 $E$ is a continuous $B$-module.
\par
We prove that $E$ is a nondegenerate $B$-module.
Let $\pi\colon B\tens E\longrightarrow E$ be the canonical map.
Let
$\Psi_{n} \in \cc$ be a sequence of positive functions such
that $\text{supp}\Psi_{n}
\longrightarrow 0$, and $\int_{G}\Psi_{n}(g)dg =1$.
Let $\Psi_{n}\otimes a$ denote the function $g \mapsto \Psi_{n}(g)a$
in $B$.
To see that $E$ is a nondegenerate $B$-module, it
suffices to show that $(\Psi_{n} \otimes a)e$ converges to $ae$
in $E$ for every $a\in A$ and $e\in E$, since then $\pi$ will have dense
image, and the null space for the action of $B$ on $E$ will
be contained in the null space for the action of $A$ on $E$.
We estimate
$$\aligned \parallel (\Psi_{n} \otimes a)e - ae \parallel_{d}
& \leq \int_{G} \Psi_{n}(g) \parallel age - ae \parallel_{d} dg\\
& \leq \sup_{g\in \text{supp}\Psi_{n}} \parallel age -ae\parallel_{d}\\
&\leq \parallel a \parallel_{m}
\sup_{g\in \text{supp}\Psi_{n}} \parallel ge -e\parallel_{k}\endaligned
\tag 5.6 $$
which tends to zero by the strong continuity of the action of $G$
on $E$.  Hence $E$ is a nondegenerate $B$-module.
\par
Now we show that $\pi \colon B\tens E \longrightarrow E$ in onto.
Since $G$ acts
differentiably on $E$, any element $e$ is a finite sum of elements
$\alpha_{f}({\tilde e}) \in E$, where $f\in \cc$
and ${\tilde e}\in E$ \cite{DM, Theorem 3.3}.  So it suffices
to show that elements of the form $\alpha_{f}({\tilde e})$
are in the image of $\pi$.
\par
Let $\tp \colon A\tens E \longrightarrow E$ be the canonical
map for the action of $A$ on $E$.
Since $\tilde \pi $ is onto by assumption,
using \cite{Tr, Theorem 45.1}
we can write
$$ {\tilde e} = \tp \biggl(\sum_{n=0}^{\infty} \lambda_{n} a_{n}
\otimes e_{n}
\biggr)  $$
where $\sum |\lambda_{n}| < 1$, and
$a_{n}\longrightarrow 0$ in $A$, $e_{n}\longrightarrow 0$ in $E$.
Then
$$\aligned
 \alpha_{f} (\te ) &= \alpha_{f} \tp \biggl(\sum_{n=0}^{\infty} \lambda_{n}
a_{n} \otimes e_{n} \biggr) \\
&= \sum_{n=0}^{\infty} \lambda_{n}  \alpha_{f} (\tp(a_n \otimes e_n )).
\endaligned
\tag 5.7 $$
Since $G$ acts differentiably on $A$, the function
 $b_n (g) = f(g)\alpha_{g} (a_n )$ is in $B$.
A simple calculation shows that $$\split
\alpha_{f} (\tp  (a_{n}\otimes e_n ))&=
\int_{G} f(g)(ga_{n}e_{n})dg \\
&=\int_{G} f(g) \alpha_{g}(a_{n})(ge_{n})dg = b_{n}e_{n}=
\pi (b_n  \otimes e_{n} ) \endsplit \tag 5.8$$
By the product rule for differentiation, and since $f\in \cc$,
we have
$$\aligned
\pa b_{n} \pa_{m, \gamma, d}&=
\int_{G} \w(g)^{m} \pa X^{\gamma} b_{n}(g) \pa_{d} dg\\
&= \int_{G} \w(g)^{m} \pa X^{\gamma}(f(g)\alpha_{g}(a_{n}))\pa_{d} dg\\
&  \leq C \sup_{|\beta|\leq |\gamma|, g\in \text{supp}(f)}
\pa X^{\beta}\alpha_{g}(a_{n}) \pa_{d} \quad{\text{$\w$ bdd on cmp sets}}\\
&  \leq C \sup_{ g\in \text{supp}(f)}
\pa \alpha_{g}(a_{n}) \pa_{k}\quad{\text{$G$ acts diff on $A$}} \\
&  \leq D
\pa a_{n} \pa_{l}\quad{\text{$G$ temp action.}} \endaligned
\tag 5.9$$
So $b_n
\longrightarrow 0$ in $B$
as $n\longrightarrow \infty$.
Hence the sum
$$ \sum_{n=0}^{\infty} \lambda_{n}
b_n \otimes e_{n}$$
converges absolutely in $B\hat \otimes E$,
and by (5.7) and (5.8) its image under
$\pi$ is $\alpha_{f}(\te )$.  We have proved that
$\pi$ is onto.  Thus $E$ is a differentiable $B$-module.
\par
{\bf Proof of Converse:\ }  We assume that $E$ is a differentiable
$B$-module.   Then $E$ is a quotient of the $B$-module $B\tens E$,
where $B$ acts on the left factor.  If we let $G$ act on $B$ by
$$(gF)(h)= \alpha_{g}(F(g^{-1}h)), \qquad g, h \in G, F\in B,\tag 5.10 $$
then the corresponding action of $G$ on $B\tens E$ on the
left factor gives rise to an action of $G$ on the quotient $E$.
Since the action (5.10) of $G$ on $B$ is both differentiable and tempered,
so is the action of $G$ on $E$.
\par
Similarly, the algebra $A$ acts on $B$ via
$$ (aF)(h) = aF(h), \qquad a \in A,\, F\in B,\, h\in G \tag 5.11$$
Using our hypothesis that $A$ is a self-differentiable Fr\'echet
algebra, we show that the action (5.11)
 makes $B$ into a differentiable $A$-module.
The action (5.10) makes
the $L^{1}$ $\w$-rapidly vanishing functions
$L_{1}^{\w}(G, A)$ \cite{Sc 1, \S 2} into an $A$-module.
By \cite{Schw, \S 5},
we may write $L_{1}^{\w}(G, A) \cong L_{1}^{\w}(G)\tens A$.
Since $A\tens L_{1}^{\w}(G, A) \cong L_{1}^{\w}(G) \tens A \tens A$,
and the map $A\tens A \longrightarrow A$ is onto, we see
that the canonical map $A\tens L_{1}^{\w}(G, A) \longrightarrow
 L_{1}^{\w}(G, A)$ is onto (the projective tensor product of
surjective maps is surjective \cite{Tr, Proposition 43.9}),
so
$L_{1}^{\w}(G, A)$
is a differentiable $A$-module.
It follows from \cite{Sc 1, Theorem A.10} that if a $G$-module $F$ is
a differentiable $A$-module, such that
the action of $G$ on $F$ commutes with
the action of $A$ on $F$, then the set of $C^{\infty}$-vectors
$\Finf$ for the action of $G$ is also a differentiable
$A$-module.  So $B$ is a differentiable $A$-module,
since it is the set of $C^{\infty}$-vectors for the
action $(gF)(h)=F(g^{-1}h)$ of $G$ on $L_{1}^{\w}(G, A)$.
\par
Since $A\tens B \longrightarrow B$ is onto,
the map $A\tens B \tens E \longrightarrow B\tens E$
is onto \cite{Tr, Proposition 43.9}, so
$B\tens E$ is a differentiable $A$-module,
and so by passing to the quotient,
we get a differentiable action of $A$ on $E$.
\par
To see the covariance of the actions of $G$ and $A$ on $E$,
it suffices to notice that (5.10) and (5.11) give covariant
actions of $G$ and $A$ on $B$.  Also, if we integrate (5.10) and
(5.11) via formula (5.4), we get $B$ acting on $B$ by left multiplication.
So the integrated form of our covariant actions of $G$ and $A$ on
$E$ will give back the original action of $B$ on $E$.
This proves the converse.
\par
Since we saw that (5.10) and (5.11) give an $\w$-tempered differentiable
covariant representation of ($G$, $A$) on $B$, we know
by the first part of the theorem that $B$ is differentiable as
a left module over itself.  To see that $B$ is self-differentiable,
it suffices to show that for every nonzero $b \in B$, $b * B \not= \bigl\{
0 \bigr\}$.    Since $A$ is nondegenerate, find $a \in A$ such
that $ba \not= 0$, and  let $\Psi_{n}\in \cc$ be as above.
Then $b* (\Psi_{n} \otimes a) \rightarrow
ba$ in $B$, so there must be some $n$ such that $b * (\Psi_{n} \otimes a)
\not= 0$.
\qed \enddemo
\proclaim{Theorem 5.12}
Let $G$ be a Lie group with sub-polynomial scale
$\w$  such that $\w_{-}$
bounds $Ad$ on $G$, and let
$A$ be a self-differentiable Fr\'echet algebra, with an
$\w$-tempered, differentiable action of $G$.
If $A$
satisfies the factorization property, then $\GLA$
is self-differentiable and satisfies the factorization property.
\endproclaim
\demo{Proof} The self-differentiability follows from
the previous theorem.
Let $E$ be a differentiable $\GLA$-module.
Then there is an associated covariant representation of ($G$, $A$)
on $E$ by the previous theorem.  Since $G$ acts differentiably on $E$, we
may apply \cite{DM, Theorem 3.3} to see that $E$ is the
finite span of $\alpha_{f}(e)$ where $f \in C_{c}^{\infty}(G)$ and
$e \in E$.  By assumption, every $e\in E$ may be written as a
finite sum of elements of the form $a {\tilde e}$ where $a \in A$ and
${\tilde e} \in E$.  Define $b(g) = f(g) \alpha_{g}(a)$.
Since $G$ acts differentiably on $A$, $b \in \GLA$.  Also,
$b {\te} = \alpha_{f}(a{\tilde e})$, so every element of
$E$ is  a finite sum of elements of the form $b{\te}$.
\qed
\enddemo
\proclaim{Corollary 5.13} Let $(M, \s, H)$ be any scaled $(G, \w)$-space,
with $\s$ a continuous, proper scale.
Then the smooth crossed product $G\rtimes^{\w} \schwmh$
is self-differentiable and satisfies the
factorization property.
\endproclaim
\demo{Proof} By Theorem 4.6, we know that $\schwmh$ satisfies the
factorization property.  Hence by Theorem 5.12, the smooth crossed
product $G\rtimes^{\w} \schwmh$ does also.
\qed
\enddemo
\par
It follows that many of the examples of smooth dense
subalgebras of transformation group C*-algebras in \cite{Sc 1, \S 5}
satisfy the factorization property.  In particular, see \cite{Sc 1, \S 5},
Examples
5.14, 5.18-20,
5.23-4, and 5.26 if $M$ is compact.
See also \cite{Sc 2, \S 10, \S 18, and Examples 3.15, 12.25}.
\enddocument